\documentclass[12pt,reqno]{amsart}
\usepackage{amsthm,amsmath}
\usepackage{enumerate}
\usepackage{fancyhdr}
\usepackage{amssymb}
\usepackage{mathrsfs}
\usepackage{verbatim}
\usepackage[body={36cc,51cc}, hmargin={2.8cm,2.8cm}]{geometry}
\usepackage{appendix}
\usepackage{color}
\usepackage[misc,geometry]{ifsym}
\usepackage{graphicx}
\usepackage{pdfsync}
 \newtheorem{theorem}{Theorem}[section]
 
 \newtheorem{lemma}[theorem]{Lemma}
 \newtheorem{proposition}[theorem]{Proposition}
 \numberwithin{equation}{section}
\allowdisplaybreaks
 \theoremstyle{definition}
 
 \theoremstyle{definition}

 \newtheorem{remark}[theorem]{Remark}
 \theoremstyle{remark}

\newcommand\norm[1]{\lVert#1\rVert}

\newcommand{\py}{\partial_y}

\newcommand{\pt}{\partial_t}
\newcommand{\p}{\partial}

\newcommand{\ga}{\gamma}
\newcommand{\les}{\lesssim}
\newcommand{\ges}{\gtrsim}
\newcommand{\lt}{\left}
\newcommand{\rt}{\right}

\begin{document}

\let\origmaketitle\maketitle
\def\maketitle{
  \begingroup
  \def\uppercasenonmath##1{} % this disables uppercasing title
  \let\MakeUppercase\relax % this disables uppercasing authors
  \origmaketitle
  \endgroup
}

\title[ Euler Equations with Damping and Solid Core]{On  Vacuum   Free Boundary Problem of the  Spherically Symmetric  Euler  Equations with  Damping and Solid Core}

\author[ Y.-L. Wang]{  Yan-Lin Wang
\\
\\
Department of Applied Mathematics,  \\  Zhejiang University of Technology \\
\Letter \  yanlin90118@163.com; yanlinwang.math@gmail.com
}
%\thanks{ \  \Letter \  yanlwang@mail.tsinghua.edu.cn ; yanlinwang.math@gmail.com}
\address[Y.-L. Wang]{Department of Applied Mathematics, Zhejiang University of Technology, Hangzhou 310023,  China.}
\email{yanlin90118@163.com;  yanlinwang.math@gmail.com}

\begin{abstract}
%In 1990s, the global well-posedness theory for vacuum  free boundary problem of the spherically symmetric compressible Navier-Stokes equations with solid core have been investigated in a series  \cite{OM1, MOM2, MOM3} by Matu\v{s}${\rm \mathring{u}}$-Ne\v{c}asov\'{a} \v{S}.,  Okada M. and Makino T..
%They proved the long-time asymptotic stability of the equilibrium with respect to small perturbation,    when $\ga>4/3$ and $A$ is sufficiently small in the state equation $p=A \rho^\ga$, with $p$ being the pressure and $\rho$ the density. As they pointed out, the assumption of  smallness  $A$  is a serious restriction. 
The global existence of smooth  solution and the long-time asymptotic stability  of the equilibrium to  vacuum free boundary problem of  the spherically symmetric  Euler equations with damping and solid core have been obtained  for arbitrary finite positive gas constant $A$ in the state equation $p=A \rho^\gamma$ with $p$ being the pressure and 
$\rho$ the density,   provided  that  $\gamma>4/3,$ initial perturbation is small  and the radius of the equilibrium $R$  is suitably  larger than the radius 
of the solid core $r_0$. Moreover,  we  obtain the  pointwise convergence from the smooth solution to the equilibrium in a  surprisingly   exponential  time-decay rate.   
The proof is mainly  based on  weighted energy method in Lagrangian coordinate.
%which is mainly  motivated by  Luo and Zeng's work [{ Comm. Pure Appl. Math. 69 (2016), no. 7, 1354–1396}] and Zeng's work [{ Arch. Ration. Mech. Anal. 226 (2017), no. 1,  33–82}] on nonlinear asymptotic stability of Barenblatt solution for vacuum  free boundary problem of Euler equations with damping.
\vspace{0.2cm}


\end{abstract}

\maketitle

\section{introduction}
The motions of the damped atmosphere surrounding a  heavy planet  is governed by the following free boundary problem
\begin{subequations}\label{3dmodel}
\begin{align}
&\pt\rho+\nabla\cdot (\rho \bold{u})=0 &&  {\rm in}\ \Omega(t),\label{3dedc-a}\\
&\rho \lt(\pt\bold{u} +\bold{u} \cdot\nabla \bold{u} \rt)+\nabla p=-\rho \bold{u}-\rho \nabla \Phi_{c}  &&  {\rm in}\ \Omega(t),\label{3dedc-b}\\
&\rho>0 && {\rm in}\ \Omega(t),\label{3dedc-c}\\
&\rho=0, && {\rm on}\ \Gamma(t)=\partial \Omega(t),\label{3dedc-d1}\\
& \bold{u}(\bar\Gamma, t)=0,\label{edc-d2}\\
&\mathcal{V}({\Gamma}(t))=\bold{u}\cdot \mathcal{N},\label{edc-e}\\
&(\rho, \bold{u})=(\rho_0, \bold{u}_0)(\bold{x})&& {\rm on}\ \Omega(0),\label{edc-f}
\end{align}
\end{subequations} 
where $(t, \bold{x})\in \mathbb{R}_+\times \mathbb{R}^3$ and $\rho, \bold{u},$ and $p$ denote the density, velocity, and pressure of gas atmosphere.  $\Omega(t)$ is the changing volume occupied by the gas at time $t,$ and $\Gamma(t)$ represents the moving vacuum boundary of the gas atmosphere. Here $\bar \Gamma, \mathcal{V}({\Gamma}(t)), \mathcal{N}$ denote, respectively,  the fixed  boundary of the solid planet, the velocity of the moving boundary $\Gamma(t),$ and the exterior unit normal vector to $\Gamma(t)$. 
 $\Phi_c$ is the gravitational potential of the   solid planet satisfying 
\begin{align}
-\Phi_c=4\pi G_0\rho_B\ast\frac{1}{|\bold{x}|}=4\pi G_0\int_B \frac{\rho_B(\bold{y})}{|\bold{x}-\bold{y}|} d\bold{y},
\end{align}
where $\rho_B$ is the nonnegative density of  the solid planet $B$ and $G_0$ is the gravitational constant. The self-gravity of the atmosphere is neglected.  That is, the total mass of the gas is relatively small compared with that of the inner solid planet. For more
information on  gaseous planets with solid core  in astrophysics, one can refer to \cite{AA, BHB, MH, MF}.  

In this paper, we assume that $\rho_B$ is spherically symmetric and   the solid planet is a  ball with radius $r_0>0.$   In spherically symmetric setting 
$$\rho(t, \bold{x})=\rho(t, r),\ \  \bold{u}(t, \bold{x})=u(t, r)\frac{\bold{x}}{r},\ \ \ {\rm with}\ \ r=|\bold{x}|\in (r_0, R(t)),$$
we rewrite the system \eqref{3dmodel} into the following 
 free boundary problem of the compressible Euler equations with damping and solid core:
\begin{subequations}\label{edc}
\begin{align}
&\rho_t+(\rho u)_r+\frac{2}{r}\rho u=0 &&  {\rm in}\ (r_0, R(t)),\label{edc-a}\\
&\rho \lt(u_t +u u_r \rt)+p_r=-\frac{ \rho g_0}{r^2}-\rho u  &&  {\rm in}\ (r_0, R(t)),\label{edc-b}\\
&\rho>0 && {\rm in}\ [r_0, R(t)),\label{edc-c}\\
&\rho(R(t),t)=0,\ u(r_0, t)=0,\label{edc-d}\\
&\dot{R}(t)=u(R(t), t),\ R(0)=R_0,\label{edc-e}\\
&(\rho, u)(r, t=0)=(\rho_0, u_0)(r)&& {\rm on}\ (r_0, R_0),\label{edc-f}
\end{align}
\end{subequations}
where  $r_0$ is the radius of the central solid core, and $g_0=G_0 M_0>0,$ with 
 $M_0=4\pi\int_0^{r_0}\rho_B s^2 ds$ being  the mass of the central core. The setup \eqref{edc-c} means that no vacuum exists between the solid core and the gas. 
In the sequel, we always assume  that 
\begin{align}p=p(\rho)= A \rho^\ga,\ \  \ga>1,\label{pressure}\end{align}
where $A$ is a finite positive constant and  $\ga$ is adiabatic exponent. 

The equilibria $(\bar\rho(r), 0)$ of  problem \eqref{edc} satisfies 
\begin{align}
 A (\bar\rho^\ga)_r=-\frac{g_0 \bar\rho}{r^2}.\label{equeqn}
\end{align}
More precisely, $\bar\rho(r)$ is given by (cf. \cite{Lin, Makino})
\begin{equation}\label{equi}
\bar\rho(r)=
\begin{cases}
\bar A\lt(\frac{1}{r}-\frac{1}{R}\rt)^{1/(\ga-1)}\ \ &r_0\leq r <R,\\
0\ \ & r\geq R,
\end{cases}
\end{equation}
where $R$ is an arbitrary number such that $R> r_0$ and  $\bar A$ is defined  by
\begin{align}\bar A=\lt(\frac{(\ga-1)g_0}{\ga A }\rt)^{1/(\ga-1)}.\label{barA}\end{align}

Note that  adiabatic exponent $\gamma$ is usually  restricted to $\ga\in(1, 2)$ for polytropic gas. For instance, $\ga\rightarrow 1_+$ is  for heavier molecules,  and $\ga=\frac{7}{5}$ for a monatomic gas and $\ga=\frac{5}{3}$ for a atomic gas. However,  for liquid (or compressed gas),  the value of   $\ga$ may be extremely  high such that $\ga \geq 2.$ (see \cite{LZ}, nonlinear stability for one dimensional  Euler equations with damping and any given $\ga>1.$ )  Indeed, the solution \eqref{equi} solves the equilibria equation \eqref{equeqn} uniquely for any given $\ga>1$ mathematically. 

The total mass $M$ of the equilibrium is given by
\begin{align}
M=4\pi \bar A \int_{r_0}^R \lt(\frac{1}{r}-\frac{1}{R}\rt)^{1/(\ga-1)} r^2 dr,
\end{align}
which is a increasing function of $R.$ It has been  clarified  in \cite{Lin, Makino} that  given  any  total mass $M\in (0,\infty)$ for $\ga\geq 4/3$, or $M\in (0, M^*)$ for $1<\ga<4/3,$  then the radius $R$ is  uniquely determined. Here
$M^*$ is given by
$$M^*=\frac{4\pi \bar A (\ga-1)}{4-3\ga}r_0^{-(4-3\ga)/(\ga-1)}.$$
Without loss of generality, we assume $A=1$ for simplicity of the symbol in the following parts.

Let $c(\rho)=\sqrt{p'(\rho)}$ be the sound speed, the condition
\begin{align}
-\infty< \p_r(c^2(\rho))<0 \ \ \ {\rm on}\ \ \ R(t)\label{phyb}
\end{align}
defines a {\it physical \ vacuum} boundary (cf. \cite{CLS, CS2, LTP, JM, LY1, LY2}),  which  is also called a vacuum boundary with physical singularity.
To capture this singularity, as in \cite{LZ, Zeng2017, Zeng2021}, we assume that  the initial density satisfies that
\begin{subequations}\label{iniden}
\begin{align}
&\rho_0(r)>0\ {\rm for}\ r_0\leq r < R_0,\ \rho(R_0)=0\  {\rm and}\ 
M=4\pi \int_{r_0}^{R_0} \rho_0(r) r^2dr,\label{ini-a}\\
&-\infty< \lt(\rho_0^{\ga-1}(r)\rt)_r<0\ {\rm at}\ r=R_0.\label{ini-b}
\end{align}
\end{subequations}
Here \eqref{phyb} and \eqref{ini-b} indicate that the sound speed is $C^{1/2}$-H\"{o}lder continuous near the vacuum boundary.

%In the following Figure 1, we  give  a profile crossing the center  of the solid ball, which is  surrounded by the atmosphere in spherical symmetry.
%\begin{figure}[!h]
%\centering
%\includegraphics[height=8cm, width=8cm]{SC202208black}
%\caption{The Profile Crossing the Center of the Solid Ball}
%\end{figure}

Clearly, the stationary solution $\bar \rho(r)$ is defined on the domain $(r_0, R)$. 
 We introduce a diffeomorphism $\eta_0: (r_0, R)\rightarrow (r_0, R_0)$ satisfying 
$$\int_{r_0}^{\eta_0(r) } r^2 \rho_0(r) dr=\int_{r_0}^r \tau^2 \bar\rho(\tau)d\tau,\ \ \ r\in (r_0, R),$$ 
which implies 
\begin{align}
\eta_0^2(r)\rho_0(\eta_0(r))\eta_{0r}(r)=r^2\bar\rho(r),\ \ \ r\in (r_0, R).\label{eta0}
\end{align}

For $y \in(r_0, R),$ we define the Lagrangian variable $\eta(y, t)$ by
\begin{align}
\eta_t(y, t)=u(\eta(y, t), t)\ \ {\rm for}\ \ t>0,\ \ {\rm and}\ \ \eta(y, 0)=\eta_0(y).\label{lvariable}
\end{align}
Note that $\eta\geq r_0$ for all $t\geq 0.$
Then the Lagrangian density and velocity are,   respectively,  defined  by
\begin{align}
\varrho(y, t)= \rho(\eta(y,t), t)\ \ {\rm and}\ \ v(y, t)=u(\eta(y,t), t), \ \ y\in(r_0, R). \label{ldv}
\end{align}
Now, by means of the Lagrangian variables defined above,  one can reformulate the system \eqref{edc} on the reference domain $(r_0, R):$
\begin{subequations}\label{ledc}
\begin{align}
&(\eta^2 \varrho)_t+\eta^2\varrho v_y /\eta_y=0 &&{\rm in}\ (r_0, R)\times(0, \infty),\label{ledc-a}\\
&\varrho v_t+(\varrho^\ga)_y/\eta_y=-\varrho v -\frac{\varrho g_0}{\eta^2} &&{\rm in}\ (r_0, R)\times(0, \infty),\label{ledc-b}\\
&v(r_0, t)=0 && {\rm on} \ \ (0, \infty),\label{ledc-c}\\
&(\varrho, v)=(\rho_0(\eta_0), u_0(\eta_0)) && {\rm on}\ \  (r_0, R)\times \{t=0\}.\label{ledc-d}
\end{align} 
\end{subequations}

It follows from \eqref{eta0} \eqref{ledc-a} that
\begin{align}
\varrho(y, t)\eta^2(y, t)\eta_y(y, t)=\varrho_0(\eta_0(y))\eta_0^2(y)\eta_{0y}(y)=y^2\bar\rho(y), \ \ y\in
(r_0, R).\label{rbr}
\end{align}
This, together with  \eqref{eta0} \eqref{lvariable} and \eqref{ldv},  enables us to  rewrite   the equations \eqref{ledc} as following  
\begin{subequations}\label{etaeqn}
\begin{align}
&\bar\rho \eta_{tt}+\bar\rho \eta_t+\lt(\frac{\eta}{y}\rt)^2 \lt[\lt(\frac{y^2\bar\rho}{\eta^2 \eta_y}\rt)^{\ga}\rt]_y+\frac{\bar\rho g_0}{\eta^2}=0 &&  {\rm in} \ \ (r_0, R)\times(0, \infty),\label{etaeqn-a}\\
&\eta(r_0, t)=r_0, \  \eta_t(r_0, t)=0, && {\rm on} \ \ (0, \infty),\label{etaeqn-b}\\
&(\eta, \eta_t)=(\eta_0, u_0(\eta_0)) && {\rm on}\ \ (r_0, R)\times \{t=0\}.\label{etaeqn-c}
\end{align}
\end{subequations}
A stationary solution to the  equation \eqref{etaeqn-a} reads
\begin{align}\bar \eta(y, t)= y, \ \ \ \bar\eta_t(y, t)=0,\ \ \ y\in (r_0, R), 
\label{beta}\end{align}
which is just  the equilibria  \eqref{equi} defined in Lagrangian version.

Let $y \zeta(y, t):=\eta(y,t)-\bar\eta(y, t)=\eta(y, t)-y.$ Then, with the aid of 
\eqref{equeqn},  we can  write the equation \eqref{etaeqn} into a perturbation 
form: 
\begin{subequations}
\begin{align}
&y \bar\rho\zeta_{tt}+y\bar\rho \zeta_t+(1+\zeta)^2\lt[\bar\rho^\ga(1+\zeta)^{-2\ga}(1+\zeta+y\zeta_y)^{-\ga}\rt]_y-\frac{ (\bar\rho^\ga)_y}{(1+\zeta)^2}=0.\label{erreqn}\\
&\zeta(r_0, t)=0, \label{zetabdr}\\
&(y\zeta, y\zeta_t)(y,0)=(\eta_0-y, u_0(\eta_0)).\label{zetaini}
\end{align}
\end{subequations}

The study of  motions of the atmosphere surrounding a planet is very important in physics and mathematics (cf. \cite{AA, BHB, MH, MF}).
In 1990s, Okada M. proved the global existence of weak solution to 1-D compressible Navier-Stokes equations with a fixed boundary and a free boundary \cite{Okada}.  Then  Matu\v{s}${\rm \mathring{u}}$-Ne\v{c}asov\'{a} \v{S}. Okada M.  and  Makino T.  proved the existence 
and uniqueness of  global  weak solution for free boundary problem of  the spherically symmetric Navier-Stokes equations with solid core in \cite{OM1, MOM2}, respectively.  Later, they \cite{MOM3} verified the  long-time   asymptotic stability ( in $L^2$ ) of  the equilibrium \eqref{equi}  with respect to small perturbation   in the regime of spherically symmetric Navier-Stokes equation with solid core in Lagrangian mass coordinate, provided that $\ga>4/3$ and the gas constant $A$ is sufficiently small (see \cite{MR} briefly).  They pointed out that  the assumption of smallness $A$ is a serious restriction but can not be removed in their scheme (cf. \cite{MOM3}). 
Then Lin S.  verified the linear stability of the equilibrium \eqref{equi} for  the spherically  symmetric Navier-Stokes-Poisson or Euler-Poisson  (Euler) with solid core when $\ga>4/3$ ($\ga>1$) but without the assumption of smallness $A$, through spectral analysis of the linearized operators.  Recently,  Makino T.  justified a long-time validity of  the approximate solution  constructed by the equilibria  plus a time-periodic solution  of the linearized equation around the equilibria   for the spherically symmetric Euler equations with solid core \cite{Makino}.  However, the nonlinear 
asymptotic stability of the equilibrium for inviscid gaseous stars  with solid core is  still an  open problem.  In this paper, we consider the global existence of strong solution to the spherically symmetric  Euler equations with damping and solid core. We will also verify the nonlinear stability of  the equilibrium \eqref{equi} under suitable assumption on the radius of the equilibrium  but without the smallness of gas constant $A.$

If there is no solid core,  some outstanding progress have been made on  the global 
existence theory and long-time  asymptotic stability for free boundary problem of  {\it viscous flow}  at present.  For instance, when $4/3<\ga<2,$  Luo T.  et al proved the long-time nonlinear  asymptotic stability of  Lane-Emden solution for the  spherically symmetric Navier-Stokes-Poisson equations with constant viscosity \cite{LXZ1, LWZ} and degenerate density-dependent viscosity \cite{LXZ2}. They applied weighted 
energy methods to overcome the extreme difficulties caused by boundary singularity or coordinate  singularity. 
  For more global existence results of viscous flow with free boundary, we  refer to \cite{HLZ, LXY, Ou, OZ, YZ, Zeng2015, Zhu} and  references therein.
  
 On the other hand, there are also some important progress on  the global existence theory and nonlinear stability for  free boundary problem of {\it inviscid flow}  without solid core. When the inviscid gas flow is    governed  by  Euler equations  with damping,  the long-time asymptotic stability of Barenblatt solution  have been obtained by Luo and Zeng in \cite{LZ}  for one dimension and by Zeng in \cite{Zeng2017} for  three dimensions in spherically symmetric setting (see also \cite{Zeng21} for almost global existence in 3-D without symmetric assumption).  
  Had\v{z}i\'c M. and Jang J \cite{HJ}, in 2019, obtained  a remarkable  result on a class of global solutions to free boundary problem of Euler-Poisson equations in three dimensions  without any  symmetric assumption for all $\ga$ in the form
  $\ga=1+\frac{1}{m}, m\in \mathbb{N} \setminus\{1\}. $ The global solutions presented in \cite{HJ} have initially small 
densities and compact supports, and they stay close to Sideris affine solutions of the Euler system. 
   More  recently,   Zeng
 H.  \cite{Zeng2021} proved the global existence for  vacuum free boundary problem of Euler equations with damping and gravity in $\mathbb{T}^n, n=1, 2, 3,$ for $\ga>1.$  Moreover,  in \cite{Zeng2021},   Zeng   obtained a surprisingly  exponential decay of  smooth solution to the stationary solution $\bar\rho(x)=(\nu(\hbar-x))^{1/(\ga-1)}$ ($\nu, \hbar$ are positive constants)  when the initial perturbation is small.  Here we give a brief comparison of the 
 background  solution mentioned above in Luo and Zeng's work.   The background solution (Barenblatt solotion) in \cite{LZ, Zeng2017} is determined by a porous media equation
$\rho_t=\Delta p(\rho),$  which is obtained through simplifying the momentum equation by the Darcy's law
$\nabla p(\rho)=-\rho u.$ The boundary of  Barenblatt solution is expanding  
sublinearly in time (cf. \cite{Barenblatt}). But the background solution for Euler equations with damping and gravity  in \cite{Zeng2021} is a stationary solution,  the boundary of which  is fixed.

However,  there are very  few results on  global existence theory and nonlinear stability  for free boundary problem of inviscid gas flow surrounding a solid core in mathematic literature.  
 At this stage, it is meaningful to consider  this issue  for  free boundary problem of   Euler equations with 
 damping and solid core in this paper.    The approach proposed in this paper can also be applied to solve the global existence and asymptotic stability  for free boundary problem of the spherically symmetric Euler-Poisson equations 
 with damping and solid core with respect to small perturbation.   The same problem for Euler (-Poisson)  equations only with solid core are  also  open problems, but they
 are extremely challenging at present, due to the much less dissipation.  If the readers  are interested
 in  the existence theory of stationary solution for  inviscid rotating stars governed by Euler-Poisson equations  with solid core,  please  refer to \cite{Wu} and  references therein.
 
 In this paper,  we will apply weighted energy method, which is shown a powerful tool to handle vacuum  free boundary problem with physical singularity (cf.\cite{JM, Jang2014, LWZ, LXZ1, LXZ2, LZ, Zeng2017, Zeng21, Zeng2021} and  references therein), to show the global-in-time existence of smooth solution and to derive the nonlinear asymptotic stability with convergence rates,  with the help of embedding of weighted Sobolev space.  The main challenges in this paper are caused by the solid core and the vacuum boundary with physical singularity. More precisely,  first, the solid core will bring a boundary effect, especially in the higher-order energy estimates, which lead to the loss of derivatives.  Hence the  approach in \cite{Zeng2021} through taking spatial derivatives directly in energy estimates  does not work smoothly   in our situation,  even though there is no coordinate singularity.  To deal with this difficulty,  we will  apply elliptic estimates  associated  with  suitable relation between the radius of the equilibrium and the radius of the solid core.   Second, for the vacuum boundary with singularity,  we  overcome this difficulty by using weighted energy estimates together with Hardy inequality,  which is mainly  inspired  by Luo and Zeng's work \cite{LZ, Zeng2017} on nonlinear asymptotic stability of Barenblatt solution for vacuum  free boundary problem of Euler equations with damping. 
In this paper,  the boundary of the  stationary solution \eqref{equi}   for the spherically  symmetric Euler equations with damping and solid core does not expand  in time.  This is at least one reason that  we can   hope a faster decay of the smooth solution   with respect  to a  small perturbation around  the equilibrium, as in \cite{Zeng2021}.

The arrangement of the remaining parts of this paper is as following: Some crucial tools for our proof and main results of  this paper will be stated in Section 2.  Next, we will prove our main results in Section 3 containing 
elliptic estimates and low-high weighted energy estimates.

%{\color{red}{Question: The global existence theory of 
 %Euler equations only with solid core but without damping seems extremely  difficult for me. 
 %Could we study the global existence of the
 % vacuum free boundary problem \eqref{edc} with damping and solid core?}}

\section{Main results}
\noindent\textbf{Notations:}
In this paper, we use $C$ to denote a generic positive  constant,  which may depend on $M, A, g_0, \ga, r_0, R,$  but independent of  the time $t$ and the initial data.  
We use $a \sim b$ to denote $C^{-1} b \leq a \leq C b,$ $a \les b$ to denote
 $a \leq C b$, and $a \ges b$ to denote $a \geq C b,$  for some generic positive
 constant $C$ defined above.  We  use $C(a)$ to represent that $C$  depends on $a$ additionally. 
 For simplicity of symbol, in the rest of the paper, we will also  use the following 
notations: 
$$\int:= \int_{r_0}^R, \ \ \norm{\cdot}:= \norm{\cdot}_{L^2(r_0, R)},\ \ 
\norm{\cdot}_{L^\infty}:=\norm{\cdot}_{L^\infty(r_0, R)}.$$
We set $\alpha=\frac{1}{\ga-1}$ and 
 $\sigma(y):=\bar\rho^{\ga-1}(y)=\bar A ^{\ga-1}(y R)^{-1}(R-y)$ for $y\in (r_0, R).$ So, $R^{-2}\bar A^{\ga-1}(R-y)\leq |\sigma(y)|\leq (r_0R)^{-1}\bar A^{\ga-1} (R-y)$ and $\sigma(y)\sim d(y)=dist (y, \p I)$.  In some occasion, we may divide the domain $I:=(r_0, R)$ into $I_l$ and $I_b$ with $I_l:=(r_0, (R+r_0)/2)$ and $I_b:=[(R+r_0)/2, R).$
 
 Near the boundary, we may use {\it Hardy inequality} (cf. \cite{KMP}) to deal with the boundary singularity.  That is, if 
 $\int_{I_b} \sigma^k(y)\lt(F^2+F_{y}^2\rt)dy<\infty,$ then,  for $k>1,$  it holds
 \begin{align}
 \int_{I_b} \sigma^{k-2}(y)F^2 dy\leq C(R, r_0,\ga,  k) \int_{I_b} \sigma^{k}(y)\lt(F^2+F_{y}^2\rt)dy. \label{hardy}
 \end{align}

We introduce the energy functionals $\mathcal{E}_j(t)$ and $\mathcal{E}_{j, i}(t)$ defined  in the following form
\begin{align}
&\mathcal{E}_{ j}(t)=\int \lt[y^4 \sigma^{\alpha}|\p_t^j(\zeta, \zeta_t)|^2+y^2\sigma^{\alpha+1}|\p_t^j(\zeta, y\zeta_y)|^2\rt](y, t) dy,\label{mej}\\
&\mathcal{E}_{j,i}(t)=\int \lt[ y^2 \sigma^{\alpha+i-1}(\pt^j\p_y^i \zeta)^2+y^4\sigma^{\alpha+i+1}(\pt^j \p_y^{i+1} \zeta)^2 \rt](y, t)dy.\label{meji}
\end{align}
The dissipation $\mathcal{D}_j(t)$ is defined by 
\begin{align}
\mathcal{D}_j(t)=\int\lt[ y^4\sigma^{\alpha} (\p_t^{j+1}\zeta)^2
+y^2\sigma^{\alpha+1}|\p_t^j(\zeta, y\zeta_y)|^2 \rt](y, t)dy.
\end{align}
Let us denote $n=4+[\alpha].$ The total energy functional  is defined by 
\begin{align}
\mathcal{E}(t):=\sum_{j=0}^n \lt( \mathcal{E}_j(t)+\sum_{i=1}^{n-j}\mathcal{E}_{j,i}(t)\rt).\label{me0815}
\end{align}
The total dissipation $\mathcal{D}(t)$ is denoted by 
\begin{align}
\mathcal{D}(t):=\sum_{j=0}^n \mathcal{D}_j(t).
\end{align}
In particular, one can derive from the verification of    Lemma 3.7  in  \cite{Zeng2017} that, for $n=4+[\alpha]$ and  $\mathcal{E}(t)$ defined in \eqref{me0815}, it holds 
\begin{align}
&\sum_{j=0}^3\norm{\pt^j\zeta(y,t)}_{L^\infty}^2+\sum_{j=0}^1\norm{\pt^j\zeta_y(y,t)}_{L^\infty}^2\notag\\
&+\sum_{0\leq j+i\leq n}\norm{y^2 \sigma^{\max\{0, \frac{2i+j-3}{2}\}}\pt^j\py^i\zeta(y, t)}_{L^\infty}^2\les \mathcal{E}(t),\ \ \ y\in(r_0, R)
\label{zeng2017}
\end{align}
provided that $\mathcal{E}(t)$ is  finite. In the derivation of \eqref{zeng2017}, the following  {\it embedding of weighted Sobolev space} play a crucial role:
\begin{align}
H^{a, b}(I)\hookrightarrow H^{b-a/2}(I),\ {\rm with}\ 
\norm{F}_{H^{b-a/2}}\leq C\norm{F}_{H^{a, b}},\label{Sembeding}
\end{align}
where the weighted Sobolev space $H^{a, b}(I)$ is given by
$$H^{a, b}(I):=\lt\{ d^{a/2}(y) F\in L^2(I): \int_I d^a(y) |\py^k F|^2dy<\infty,\ \ 0\leq k\leq b \rt\},$$
with  the norm $\norm{F}^2_{H^{a,b}(I)}=\sum_{k=0}^b \int_I d^a(y) |\py^k F|^2dy.$

\noindent\textbf{A Priori Assumption}.
Let $\eta(y,t)=y\zeta(y)+y$ be a solution to problem \eqref{etaeqn} satisfying the {\it a priori assumption}:
\begin{align}
\mathcal{E}(t)\leq \epsilon^2_0,\ \ t\in[0, T],\label{assumption}
\end{align}
for some suitably small fixed positive constant $\epsilon_0$ independent  of $t.$ 
Therefore, it follows from \eqref{zeng2017} and \eqref{assumption} that  for $y\in (r_0, R), t\in[0, T],$
\begin{align}
&\sum_{j=0}^3\norm{\pt^j\zeta(y,t)}_{L^\infty}^2+\sum_{j=0}^1\norm{\pt^j\zeta_y(y,t)}_{L^\infty}^2\notag\\
&+\sum_{0\leq j+i\leq 4+[\alpha]}\norm{y^2 \sigma^{\max\{0, \frac{2i+j-3}{2}\}}\pt^j\py^i\zeta(\cdot,t)}_{L^\infty}^2\les \mathcal{E}(t)\leq\epsilon_0^2.\label{apriori}
\end{align}

Notice that the  coordinate singularity will not appear in this paper, due to \eqref{edc-c} and  the inner solid core (cf. \cite{M1993, MOM3}). 
Instead, we should pay attention to the boundary effect contacting  the solid sphere.  Now we state the main results as following.
\begin{theorem}\label{thm1}
Suppose $\ga>\frac{4}{3}$ in the barotropic gas pressure \eqref{pressure} and $r_0< R\leq \frac{4}{3-\alpha}r_0<\infty$ with $\alpha=\frac{1}{\ga-1},$  where $R$  is the finite radius of the equilibria and $r_0$ is the radius of the solid core.  There exists a positive constant $\epsilon\in(0, \epsilon_0]$  such that if $\mathcal{E}(0)\leq \epsilon^2,$ then problem \eqref{etaeqn} admits a global smooth solution in $(r_0, R)\times[0,\infty)$ satisfying 
\begin{align}
\mathcal{E}(t)\les e^{-\delta t} \mathcal{E}(0),\ \ \ t\geq 0,\label{thm1-1}
\end{align}
for some suitably small positive constant  $\delta$ independent of $t.$
\end{theorem}

\begin{remark}
As explained in last section, the radius $R$ is uniquely determined by  the total mass $M=4\pi \bar A\int_{r_0}^R \lt(\frac{1}{r}-\frac{1}{R}\rt)^{1/(\ga-1)}r^2 dr$ of the atmosphere surrounding planet.  Hence, there always exists some given  finite  total mass  $M>0$ such that  $r_0 <R \leq  \frac{4}{3-\alpha}r_0<  \infty$ holds.  Here the upper bound of $R$ is merely  a technical assumption in our paper, due to the jump discontinuity of the density $\rho(t, r)$ contacting the solid sphere (see elliptic estimates in section \ref{ellip3.1}). When $\ga>\frac{4}{3},$  it  believes that Theorem \ref{thm1} holds for any finite total mass $0<M<+\infty,$ i.e.,  $r_0<R<+\infty.$ Usually, we restrict $\ga\in(1, 2)$ for isentropic  ploytropic gas, but Theorem \ref{thm1} also holds for any  given $\ga\geq 2$ mathematically. 
\end{remark}

As a corollary of Theorem \ref{thm1}, we have the following convergence results of solutions to the original vacuum 
free boundary problem \eqref{edc} concerning the vacuum free boundary $R(t),$ density $\rho$ and  velocity $u.$

\begin{theorem}\label{thm2}
Suppose $\ga>\frac{4}{3}$ and $r_0<  R \leq  \frac{4}{3-\alpha}r_0<\infty.$ 
 There exists a positive constant $\epsilon\in(0, \epsilon_0]$  such that if $\mathcal{E}(0)\leq \epsilon^2,$ then problem \eqref{edc} admits a global smooth solution $(\rho, u, R(t))$ for $t\in [0, \infty)$ satisfying 
 \begin{align}
 &\lt|\rho(t, \eta(t,y))-\bar\rho(y)\rt| \les \bar A
 \lt(\frac{1}{y}-\frac{1}{R}\rt)^{1/(\ga-1)}\sqrt{e^{-\delta t}\mathcal{E}(0)},\label{rhodec}\\
 &\lt|u(t, \eta(t, y))\rt|+\lt|R(t)-R\rt|+\sum_{1 \leq \iota\leq 3}\lt|\frac{d^\iota R(t)}{dt^\iota}\rt|\les \sqrt{e^{-\delta t}\mathcal{E}(0)},\label{urdec}
 \end{align}
 for all $y\in(r_0, R) $ and some suitably small  positive constant $\delta$ independent of $t,$  where $\bar A$ is a positive  constant given in \eqref{barA}.
\end{theorem}

%\begin{remark}
If the self-gravity of the atmosphere surrounding planet can not  be neglected,   then  the spherically symmetric  motion of the atmosphere  is governed by the following    system:
\begin{subequations}\label{epc}
\begin{align}
&\rho_t+(\rho u)_r+\frac{2}{r}\rho u=0 &&  {\rm in}\ (r_0, R(t)),\label{epc-a}\\
&\rho \lt(u_t +u u_r \rt)+p_r=-\rho u -\frac{\rho}{r^2}\lt(g_0+4\pi G\int_{r_0}^r \rho(t,\tau)\tau^2 d\tau\rt) &&  {\rm in}\ (r_0, R(t)),\label{epc-b}\\
&\rho>0 && {\rm in}\ [r_0, R(t)),\label{epc-c}\\
&\rho(R(t),t)=0,\ u(r_0, t)=0,\label{epc-d}\\
&\dot{R}(t)=u(R(t), t),\ R(0)=R_0,\label{epc-e}\\
&(\rho, u)(r, t=0)=(\rho_0, u_0)(r)&& {\rm on}\ (r_0, R_0),\label{epc-f}
\end{align}
\end{subequations}
which is  Euler-Poisson system with damping and solid core.
% the momentum equation \eqref{edc-b} changes into
%\begin{align}\rho (u_t+uu_r)+p_r=-\rho u-\frac{\rho}{r^2}\lt(g_0+4\pi G\int_{r_0}^r \rho(t,\tau)\tau^2 d\tau\rt),\label{epc}\end{align}
Here $G$ is the gravitational constant  for the atmosphere  surrounding planet. 

The equilibrium $(\bar\rho_*, 0)$ of  \eqref{epc} with $g_0>0$
satisfies
\begin{align}
A(\bar\rho_*^\ga)_r=-\frac{\bar\rho_*}{r^2}\lt(g_0+4\pi G\int_{r_0}^r \bar\rho_* (t,\tau)\tau^2 d\tau\rt),\label{gravity1}
\end{align}
which has a unique solution  if $\ga\geq4/3$, $g_0>0$ and the finite  total mass of gas $M'=4\pi \int_{r_0}^{R_G} \bar\rho_* r^2dr>0$
(cf. \cite{KL}). Here $R_G$ is the first zero of $\bar\rho_*(r).$ We assume $g_0>0$ for the case with solid core in this paper. However,  if there is no solid core in \eqref{gravity1}, i.e. 
$g_0=0,$  then  \eqref{gravity1} is  reduced  to Lane-Emden equation, which also  has a  unique and compact supported 
 solution provided  that $ \ga\in (4/3,2)$ and $0<M<+\infty$ (cf. \cite{Chandr,  LY87, Lin, LS, Rein}). 
 
For system \eqref{epc}, similarly, we define the Lagrangian variable $\eta(t, y)$ on the  reference domain $y\in (r_0, R_G)$ and obtain  the following theorem.
\begin{theorem}\label{thm3}
Suppose $\ga>\frac{4}{3}$ and $r_0<  R_G \leq  \frac{4}{3-\alpha}r_0<\infty.$ 
 There exists a small spherically symmetric  initial perturbation around the equilibrium $(\bar \rho_*, 0)$ defined in \eqref{gravity1}  such that the problem \eqref{epc} admits a global smooth solution $(\rho, u, R(t))$ for $t\in [0, \infty)$ satisfying 
 \begin{align}
 &\lt|\rho(t, \eta(t,y))-\bar\rho_*(y)\rt| \les 
 \lt(y-R_G\rt)^{1/(\ga-1)}\sqrt{e^{-\delta t}\mathcal{E}(0)},\label{rhodec}\\
 &\lt|u(t, \eta(t, y))\rt|+\lt|R(t)-R_G\rt|+\sum_{1 \leq \iota\leq 3}\lt|\frac{d^\iota R(t)}{dt^\iota}\rt|\les \sqrt{e^{-\delta t}\mathcal{E}(0)},\label{urdec}
 \end{align}
 for all $y\in(r_0, R_G) $ and some suitably small  positive constant $\delta$ independent of $t,$  where $\mathcal{E}(0)$ represents small initial perturbation around $(\bar\rho_*, 0)$ defined in the form \eqref{me0815}.
\end{theorem}

%For every 
%$M'=4\pi\int_{r_0}^{R_{G}}\bar\rho_* r^2dr<\infty$, where $R_G$ is the radius of the equilibrium \eqref{gravity1},
%if $\ga>4/3$ and  $ r_0< R_G \leq \frac{4}{3-\alpha}r_0<\infty$, the global existence theory and long-time  asymptotic stability of the equilibrium  for the spherically symmetric Euler-Poisson equations with damping and solid core  can also be  obtained with respect to small perturbation,   by a similar argument.  
%\end{remark}

\begin{remark}
We  do not need  the condition of  smallness $A$ assumed  in \cite{MOM3} (which is for viscous case)  in our setup.  
In \cite{MOM3},  Matu\v{s}${\rm \mathring{u}}$-Ne\v{c}asov\'{a} \v{S}.,  Okada M. and Makino T.  conjectured that the assumption $\ga>4/3$ can be removed 
for the  asymptotic stability of the 
equilibrium  for  the spherically 
symmetric Navier-Stokes equations with solid core when the self-gravitation of gas  is neglected.  Similarly, 
relaxing  the  assumption  $\ga>4/3$ in Theorem \ref{thm1} to $\ga>1$  is also an unsolved  problem. But, $\ga>4/3$ may be essential for the stability of the equilibrium of Euler-Poisson equations 
with damping and solid core \eqref{epc} in Theorem \ref{thm3} (cf.\cite{KL, MOM3}).

\end{remark}

\section{Proof of Theorem \ref{thm1} -- Theorem \ref{thm3}}
The proof of Theorem \ref{thm1} and Theorem \ref{thm2} consists of elliptic estimates (controlling 
   spatial derivatives by temporal derivatives)  and  weighted energy estimates (dealing with 
    temporal derivatives).  Indeed,  it is shown this strategy  works
    in dealing with the boundary effect near the solid core and  singularity near the  vacuum boundary  in this paper.  Theorem \ref{thm3} can be proved similarly.  The smallness of $A,$  which is assumed in  nonlinear asymptotic stability of the equilibrium for the  spherically symmetric  motion of viscous gas surrounding a solid sphere (cf. \cite{MOM3}),   is not required in this paper if the radius of the equilibrium is suitably larger than the radius of the inner solid ball.

\subsection{Elliptic Estimates}\label{ellip3.1}

To begin with, we state the main result of this subsection as following.
\begin{proposition}\label{pro1}
Suppose that \eqref{apriori} holds for a suitably small constant $\epsilon_0\in (0, 1).$  If $\ga>4/3$ and 
$r_0< R \leq \frac{4}{3-\alpha}r_0<\infty,$ then it holds that 
\begin{align}
\mathcal{E}_{j, i}(t)\les \sum_{\iota=0}^{i+j} \mathcal{E}_\iota(t),\ \ \ \forall t\in [0, T],\label{proest}
\end{align}
when $j\geq 0, i\geq 0, i+j\leq n.$ Moreover, the elliptic estimates  \eqref{proest} also hold for any  $\ga\in (1, 4/3]$  and $r_0<R<\infty.$
\end{proposition}

The proof of the elliptic estimates in Proposition \ref{pro1} consists of Lemma \ref{le-elliptic1} and Lemma \ref{le2}.

Note that $\sigma_y(y)=-{\bar A}^{\ga-1}/y^2, y\in (r_0, R).$ 
Dividing equation \eqref{erreqn} by $\bar \rho,$ we obtain
\begin{align}
&y\zeta_{tt}+y\zeta_t+\sigma(1+\zeta)^2\lt[(1+\zeta)^{-2\ga}(1+\zeta+y\zeta_y)^{-\ga}\rt]_y\notag\\
&\ +\frac{\ga}{\ga-1}\sigma_y\lt[(1+\zeta)^{2-\ga}(1+\zeta+y\zeta_y)^{-\ga}-(1+\zeta)^{-2}\rt]=0.\label{dirho}
\end{align}
Here we have
\begin{align}
&(1+\zeta)^2\lt[(1+\zeta)^{-2\ga}(1+\zeta+y\zeta_y)^{-\ga}\rt]_y\notag\\
&=-\ga (4\zeta_y+y\zeta_{yy})+J_1,\notag\\
&(1+\zeta)^{2-\ga}(1+\zeta+y\zeta_y)^{-\ga}-(1+\zeta)^{-2}\notag\\
&=-\ga y\zeta_y+(4-3\ga)\zeta+J_2,\notag
\end{align}
where $J_1, J_2$ are given by
\begin{align}
J_1:=&-2\ga \zeta_y\lt[(1+\zeta)^{1-2\ga}(1+\zeta+y\zeta_y)^{-\ga}-1\rt]\notag\\
&-\ga (2\zeta_y+y\zeta_{yy})\lt[(1+\zeta)^{2-2\ga}(1+\zeta+y\zeta_y)^{-\ga-1}-1\rt],\label{j1}\\
J_2:=&(1+\zeta)^{2-2\ga}(1+\zeta+y\zeta_y)^{-\ga}-(1+\zeta)^{-2}\notag\\
&+\ga y\zeta_y-(4-3\ga)\zeta.\label{j2}
\end{align}
This, together with the Taylor expansion and the smallness of $\zeta$ and $y\zeta_y$, implies that 
\begin{align}
&|J_1|\les \epsilon_0(|y\zeta_{yy}|+|\zeta_y|),\label{j1est}\\
&|J_2|\les \epsilon_0(|y\zeta_y|+|\zeta|).\label{j2est}
\end{align}
Then we can write \eqref{dirho} as following
\begin{align}
\ga\lt[y\sigma \zeta_{yy}+4\sigma \zeta_y+\frac{\ga}{\ga-1}y \sigma_y \zeta_y\rt]=y\zeta_{tt}+y\zeta_t+\frac{\ga(4-3\ga)}{\ga-1}\sigma_y\zeta+\sigma J_1+\frac{\ga}{\ga-1}\sigma_y J_2.\label{reelliptic}
\end{align}

\begin{lemma}[Lower-Order Elliptic Estimates]\label{le-elliptic1}
Suppose that \eqref{apriori} holds for a suitably small constant $\epsilon_0\in (0, 1).$ If $\ga>4/3$ and 
$r_0< R\leq \frac{4}{3-\alpha}r_0=(1+\frac{\ga}{3\ga-4})r_0<\infty$,  then it holds that 
\begin{align}
\mathcal{E}_{0, 1}(t)\les \sum_{\iota=0}^{1} \mathcal{E}_\iota(t),\ \ \ \forall t\in [0, T]. \label{me01}
\end{align}
\end{lemma}
\begin{proof}
Square  the product of  \eqref{reelliptic} and $y \sigma^{\alpha/2}$ and then integrate the resulting formula with respect to spatial variable to obtain
\begin{align}
&\lt\| y^2 \sigma^{1+\frac{\alpha}{2}}\zeta_{yy}+4 y \sigma^{1+\frac{\alpha}{2}}\zeta_y+(1+\alpha)y^2 \sigma^{\frac{\alpha}{2}}\sigma_y \zeta_y\rt\|^2\notag\\
&\les \mathcal{E}_1(t)+\lt\| y\sigma^{1+\frac{\alpha}{2}} J_1\rt\|^2
+\lt\| y\sigma^{\frac{\alpha}{2}} \sigma_yJ_2\rt\|^2+\frac{\bar A^{2(\ga-1)}}{r_0^4}\lt\|y \sigma^{\frac{\alpha}{2}} \zeta\rt\|^2\notag\\
&\les  \mathcal{E}_1(t) +\epsilon^2_0\lt(\lt\|y^2\sigma^{1+\frac{\alpha}{2}}\zeta_{yy}\rt\|^2+
\lt\|y \sigma^{1+\frac{\alpha}{2}}\zeta_{y}\rt\|^2+\lt\|y^2 \sigma^{\frac{\alpha}{2}}\sigma_y\zeta_{y}\rt\|^2
\rt)\notag\\
&\ \ \ +\frac{\bar A^{2(\ga-1)}}{r_0^4} (1+\epsilon_0^2)\lt\|y \sigma^{\frac{\alpha}{2}} \zeta\rt\|^2,\label{le1-1}
\end{align}
where the definition of $\mathcal{E}_1$ and the fact
 \begin{align}\bar A^{\ga-1}/R^2 \leq |\sigma_y|\leq \bar A^{\ga-1}/r_0^2\label{sgy}\end{align}
 have been used for the
 first inequality, and \eqref{j1est} \eqref{j2est} for the second inequality.  The last term in \eqref{le1-1} satisfies
 \begin{align}
 \lt\|y \sigma^{\frac{\alpha}{2}} \zeta\rt\|^2
 &=\int_{r_0}^{(R+r_0)/2}y^2 \sigma^\alpha \zeta^2 dy+\int_{(R+r_0)/2}^{R}y^2 \sigma^\alpha \zeta^2 dy\notag\\
 &\les \int_{r_0}^{(R+r_0)/2}y^2 \sigma^{1+\alpha} \zeta^2 dy+\int_{(R+r_0)/2}^{R}y^4 \sigma^\alpha \zeta^2 dy
 \les \mathcal{E}_0(t).\label{le1-2}
 \end{align}
 Then it follows from \eqref{le1-1} and \eqref{le1-2} that 
 \begin{align}
&\lt\| y^2 \sigma^{1+\frac{\alpha}{2}}\zeta_{yy}+4 y \sigma^{1+\frac{\alpha}{2}}\zeta_y+(1+\alpha)y^2 \sigma^{\frac{\alpha}{2}}\sigma_y \zeta_y\rt\|^2\notag\\
&\les \frac{\bar A^{2(\ga-1)}}{r_0^4} \mathcal{E}_0(t) + \mathcal{E}_1(t) +\epsilon^2_0\lt(\lt\|y^2\sigma^{1+\frac{\alpha}{2}}\zeta_{yy}\rt\|^2+
\lt\|y \sigma^{1+\frac{\alpha}{2}}\zeta_{y}\rt\|^2+\lt\|y^2 \sigma^{\frac{\alpha}{2}}\sigma_y\zeta_{y}\rt\|^2
\rt).\label{le1-3}
\end{align}

Now we analyze the left-hand side of \eqref{le1-3}, which reads
\begin{align}
&\lt\| y^2 \sigma^{1+\frac{\alpha}{2}}\zeta_{yy}+4 y \sigma^{1+\frac{\alpha}{2}}\zeta_y+(1+\alpha)y^2 \sigma^{\frac{\alpha}{2}}\sigma_y \zeta_y\rt\|^2\notag\\
&=\norm{y^2 \sigma^{1+\frac{\alpha}{2}}\zeta_{yy}}^2+16\norm{y \sigma^{1+\frac{\alpha}{2}}\zeta_{y}}^2
+(1+\alpha)^2\norm{y^2 \sigma^\frac{\alpha}{2}\sigma_y \zeta_y}^2\notag\\
&\ \ \ +\int \lt[4y^3 \sigma^{2+\alpha}+(1+\alpha)y^4 \sigma^{1+\alpha} \sigma_y\rt](\zeta_y^2)_y dy
+8(1+\alpha)\int y^3 \sigma^{1+\alpha} \sigma_y \zeta_y^2 dy.\label{left}
\end{align}
By virtue of $\sigma_y(y)=-\bar A^{\ga-1}/y^2, y\in(r_0, R)$, for the last term in \eqref{left},  we obtain
\begin{align}
8(1+\alpha)\int y^3 \sigma^{1+\alpha} \sigma_y \zeta_y^2 dy\les& \frac{\bar A^{\ga-1}}{r^3_0}\int y^4 \sigma^{1+\alpha}  \zeta_y^2 dy\les \frac{\bar A^{\ga-1}}{r^3_0} \mathcal{E}_0(t).
\end{align}
For the last second term in \eqref{left}, with the aid of $\bar\rho(R)=0$,  integration by parts gives that
\begin{align}
\int &\lt[4y^3 \sigma^{2+\alpha}+(1+\alpha)y^4 \sigma^{1+\alpha} \sigma_y\rt](\zeta_y^2)_y dy\notag\\
=&-\{[4y^3 \sigma^{2+\alpha}+(1+\alpha)y^4 \sigma^{1+\alpha} \sigma_y]\zeta_y^2\}|_{y=r_0} \notag\\
&-\int \lt[4y^3 \sigma^{2+\alpha}+(1+\alpha)y^4 \sigma^{1+\alpha} \sigma_y\rt]_y \zeta_y^2 dy\notag\\
=&:L_1+L_2.
\end{align}
If $$\ga>\frac{4}{3}\ \ ({\rm i.e.\ \alpha<3})\ \ {\rm and}\ \ r_0< R\leq \frac{4}{3-\alpha}r_0= (1+\frac{\ga}{3\ga-4})r_0,$$
then
\begin{align}
L_1=-\sigma^{1+\alpha}(r_0)r_0^3\bar A^{\ga-1}\lt[\frac{3-\alpha}{r_0}-\frac{4}{R}\rt]\lt(\zeta_y(r_0)\rt)^2\geq 0. \label{L1}
\end{align}
(When $\alpha\geq 3,$ i.e. $1<\ga\leq 4/3,$ it also holds $L_1> 0$. Here we only consider $\ga>4/3$, due to $\ga>4/3$ is essential for our energy estimates in next section.)
Moreover, it holds 
\begin{align}
L_2\geq -12 \norm{y \sigma^{1+\frac{\alpha}{2}}\zeta_{y}}^2-(1+\alpha)^2\norm{y^2 \sigma^\frac{\alpha}{2}\sigma_y \zeta_y}^2- C\lt(\bar A^{\ga-1}, r_0^{-1}\rt) \mathcal{E}_0(t).\label{L2}
\end{align}

Then it follows from \eqref{le1-3}-\eqref{L2} that 
\begin{align}
&\norm{y^2 \sigma^{1+\frac{\alpha}{2}}\zeta_{yy}}^2+4\norm{y \sigma^{1+\frac{\alpha}{2}}\zeta_{y}}^2\notag\\
&\les \mathcal{E}_0+\mathcal{E}_1+\epsilon^2_0\lt(\lt\|y^2\sigma^{1+\frac{\alpha}{2}}\zeta_{yy}\rt\|^2+
\lt\|y \sigma^{1+\frac{\alpha}{2}}\zeta_{y}\rt\|^2+\lt\|y^2 \sigma^{\frac{\alpha}{2}}\sigma_y\zeta_{y}\rt\|^2
\rt),\label{le1-4}
\end{align}
which, together with \eqref{le1-3}, implies
\begin{align}
&\norm{y^2 \sigma^{1+\frac{\alpha}{2}}\zeta_{yy}}^2+\norm{y \sigma^{1+\frac{\alpha}{2}}\zeta_{y}}^2+\norm
{y^2 \sigma^{\frac{\alpha}{2}}\sigma_y \zeta_y}^2\notag\\
&\les \mathcal{E}_0+\mathcal{E}_1+\epsilon^2_0\lt(\lt\|y^2\sigma^{1+\frac{\alpha}{2}}\zeta_{yy}\rt\|^2+
\lt\|y \sigma^{1+\frac{\alpha}{2}}\zeta_{y}\rt\|^2+\lt\|y^2 \sigma^{\frac{\alpha}{2}}\sigma_y\zeta_{y}\rt\|^2
\rt).\label{le1-5}
\end{align}
This, with the aid of the smallness of $\epsilon_0$, implies 
\begin{align}
&\norm{y^2 \sigma^{1+\frac{\alpha}{2}}\zeta_{yy}}^2+\norm{y \sigma^{1+\frac{\alpha}{2}}\zeta_{y}}^2+\norm
{y^2 \sigma^{\frac{\alpha}{2}}\sigma_y \zeta_y}^2\les \mathcal{E}_0+\mathcal{E}_1.\label{le1-6}
\end{align}
Due to $\sigma_y=-\bar A^{\ga-1}/y^2, y\in(r_0,R),$ we easily obtain 
\begin{align}
\norm{y \sigma^{\frac{\alpha}{2}}\zeta_y}^2\leq R^2\norm{ \sigma^{\frac{\alpha}{2}}\zeta_y}^2\leq  \frac{R^2}{\bar A^{2\ga-2}}\norm
{y^2 \sigma^{\frac{\alpha}{2}}\sigma_y \zeta_y}^2\les \mathcal{E}_0+\mathcal{E}_1.\label{le1-7}
\end{align}
Therefore, \eqref{me01} follows from \eqref{le1-6} and \eqref{le1-7} and  the definition of $\mathcal{E}_{0,1}.$ The proof of this lemma is completed. 

 \end{proof}

%{\color{red}{How to deal with the boundary terms ? NONNEGATIVE!  The smallness of  gas constant $A$  does not need }}

Next, we process the higher-order elliptic estimates. Let us denote $C_m^j$ as the binomial coefficients:
 $$C_m^j=\frac{m!}{j!(m-j)!}, \ \ \ 0\leq j\leq m.$$
For $i\geq 1$ and $j\geq 0,$ we act $\pt^j\p_y^{i-1}$ on 
\eqref{reelliptic} to obtain that 
\begin{align}
&\ga \lt[y \sigma \pt^j\p_y^{i+1}\zeta+(i+3)\sigma \pt^j \p_y^i \zeta+ (\alpha+i)y \sigma_y \pt^j\py^i \zeta\rt]\notag\\
&=y \pt^{j+2}\py^{i-1}\zeta+y\pt^{j+1}\py^{i-1}\zeta +Q_1+Q_2,\label{highell}
\end{align}
where 
\begin{align}
{Q}_1:=&-\ga\pt^j\Big[\sum_{\iota=2}^{i-1}C_{i-1}^\iota \py^{\iota}(y\sigma)\py^{i+1-\iota}\zeta +4 \sum_{\iota=1}^{i-1}C_{i-1}^{\iota}\py^\iota\sigma \py^{i-\iota}\zeta\notag\\
&+(1+\alpha)\sum_{\iota=1}^{i-1}C_{i-1}^{\iota}\py^\iota(y\sigma)\py^{i-\iota}\zeta \Big]+(i-1)\pt^{j+2}\py^{i-2}\zeta\notag\\
&+(i-1)\pt^{j+1}\py^{i-2}\zeta+(1+\alpha)(4-3\ga)\sum_{\iota=1}^{i-1}C_{i-1}^{\iota}\py^{\iota+1}\sigma \pt^j\py^{i-\iota}\zeta,\label{q1}\\
{Q}_2:=&\py^{i-1}(\sigma \pt^j J_1)+(1+\alpha)\py^{i-1}(\sigma_y \pt^j J_2).\label{q2}
\end{align}
Here $J_1, J_2$ are given in \eqref{j1} \eqref{j2}.  Note that $\sum_{\iota=1}^{i-1}$ and $\sum_{\iota=2}^{i-1}$ indicate to be zero when $i=1$ and $i=1, 2,$ respectively. 

Multiply the above equation \eqref{highell} by $y \sigma^{(\alpha+i-1)/2},$ and then integrate  the quadratic form of the   product  with respect to spacial variable to obtain
\begin{align}
&\norm{y^2 \sigma^{\frac{\alpha+i+1}{2}}\pt^j\py^{i+1}\zeta+(i+3)y\sigma^{\frac{\alpha+i+1}{2}}\pt^j \py^i \zeta
+(\alpha+i)y^2 \sigma^{\frac{\alpha+i-1}{2}}\sigma_y \pt^j\py^i\zeta}^2\notag\\
&\les \norm{y^2 \sigma^{\frac{\alpha+i-1}{2}}\pt^{j+2}\py^{i-1}\zeta}^2+\norm{y^2 \sigma^{\frac{\alpha+i-1}{2}}\pt^{j+1}\py^{i-1}\zeta}^2\notag\\
&\ \ \ \ +\norm{y\sigma^{\frac{\alpha+i-1}{2}}Q_1}^2+\norm{y \sigma^{\frac{\alpha+i-1}{2}}Q_2}^2.\label{squ}
\end{align}
Similar to the derivation of \eqref{left}-\eqref{le1-4}, we obtain that the left-hand
side of \eqref{squ} satisfies 	
\begin{align}
&\norm{y^2 \sigma^{\frac{\alpha+i+1}{2}}\pt^j\py^{i+1}\zeta+(i+3)y\sigma^{\frac{\alpha+i+1}{2}}\pt^j \py^i \zeta
+(\alpha+i)y^2 \sigma^{\frac{\alpha+i-1}{2}}\sigma_y \pt^j\py^i\zeta}^2\notag\\
&\ges \mathcal{E}_{j,i}(t)-C \norm{y^2 \sigma^{\frac{\alpha+i}{2}}\pt^j\py^{i}\zeta}^2,\label{Oct25}\\
&{\rm if}\ \ \ga>4/3 \ \ {\rm and}\ \  R\leq \frac{3+i}{3-\alpha}r_0=\lt(1+\frac{i(\ga-1)+1}{3\ga-4}\rt)r_0.\notag
\end{align}
Indeed, \eqref{Oct25} holds  immediately  for  any $\ga\in(1, 4/3]$ and $r_0<R<\infty.$
 Then, for  $r_0< R\leq \frac{3+i}{3-\alpha}r_0, i=1, 2, \cdots, n,$ and $\ga>4/3$,  we have 
 \begin{align}
 \mathcal{E}_{j,i}(t)\les& \norm{y^2 \sigma^{\frac{\alpha+i}{2}}\pt^j\py^{i}\zeta}^2+
 \norm{y^2 \sigma^{\frac{\alpha+i-1}{2}}\pt^{j+2}\py^{i-1}\zeta}^2\notag\\
 &+\norm{y^2 \sigma^{\frac{\alpha+i-1}{2}}\pt^{j+1}\py^{i-1}\zeta}^2
 +\norm{y\sigma^{\frac{\alpha+i-1}{2}}Q_1}^2+\norm{y \sigma^{\frac{\alpha+i-1}{2}}Q_2}^2. \label{ejiq12}
 \end{align}
 
 To end the proof of the Proposition \ref{pro1}, we will prove the following higher-order elliptic estimates.
 \begin{lemma}[Higher-Order Elliptic Estimates]\label{le2}
 Suppose that \eqref{apriori} holds for a suitably small constant $\epsilon_0\in (0, 1).$ If $\ga>4/3$ and 
$r_0<R\leq \frac{4}{3-\alpha}r_0<\infty$,  then, for $j\geq 0, i\geq 0$ and $1\leq i+j\leq n$,  it holds that 
\begin{align}
\mathcal{E}_{j, i}(t)\les \sum_{\iota=0}^{i+j} \mathcal{E}_\iota(t),\ \ \ \forall t\in [0, T]. \label{ellmeji}
\end{align}

 \end{lemma}

\begin{proof}
Inspired by the elliptic estimates in \cite{LZ, Zeng2017},  we plan to prove this lemma by mathematical induction. 
First, we notice that 
\begin{align}
\mathcal{E}_{j,0}(t)=&
\int \lt[ y^2 \sigma^{\alpha-1}(\pt^j \zeta)^2+y^4\sigma^{\alpha+1}(\pt^j \p_y \zeta)^2 \rt](y, t)dy\notag\\
\les &\int_{r_0}^{(r_0+R)/2} y^4\sigma^\alpha (\pt^j \zeta)^2dy +\int_{(r_0+R)/2}^R y^4\sigma^{\alpha+1}(\pt^j\zeta_y)^2dy+\mathcal{E}_j(t)\notag\\
\les & \mathcal{E}_j(t),\ \ \ j\geq 0.\label{mej0}
\end{align}

It has been verified in Lemma \ref{le-elliptic1}
 that \eqref{ellmeji} holds for $i+j=1.$ Then, for $1\leq k \leq n-1,$ we  make the induction hypothesis that 
 \eqref{ellmeji} holds for all  $j\geq 0, i\geq 0$ and $i+j\leq k,$ i.e.
 \begin{align}
\mathcal{E}_{j, i}(t)\les \sum_{\iota=0}^{i+j} \mathcal{E}_\iota(t),\ \ \ j\geq0,\ \  i\geq 0,\ \  i+j\leq k.\label{hypo}
\end{align}
Then it suffices to prove \eqref{ellmeji} for $j\geq 0, i\geq 0$ and $j+i=k+1.$ That is, we only need to  prove
\begin{align*}
\mathcal{E}_{k+1-\iota, \iota}(t)\les \sum_{\iota=0}^{k+1} \mathcal{E}_\iota(t)
\end{align*}
from $\iota=0$ to $\iota=k+1$ step by step. 

In what follows, we naturally assume $j\geq 0, i\geq 1$ and $i+j=k+1\leq n$ ( for the case $i=0,$ it has been shown in \eqref{mej0}).  In view of the last two terms in \eqref{ejiq12}, we have to estimate $Q_1$ and $Q_2$ to achieve our goal. For $Q_1,$ when $i=1,$  we have $Q_1=0.$ When $i\geq 2,$  it follows from \eqref{sgy} that
\begin{align}
|Q_1| \les & \sum_{\iota=1}^{i-1}|\pt^j\py^\iota \zeta|+(i-1)(|\pt^{j+2}\py^{i-2}\zeta|+|\pt^{j+1}\py^{i-2}|).\notag
\end{align}
This implies that
\begin{align}
\norm{y \sigma^{\frac{\alpha+i-1}{2}}Q_1}^2\les \sum_{\iota=1}^{i-1}\norm{y \sigma^{\frac{\alpha+i-1}{2}}\pt^j\py^\iota \zeta}^2+(i-1)^2\sum_{\iota=j+1}^{j+2}\norm{y \sigma^{\frac{\alpha+i-1}{2}}\pt^\iota\py^{i-2}}^2.
\end{align}
Hence we have 
\begin{align}
\norm{y \sigma^{\frac{\alpha+i-1}{2}}Q_1}^2\les \sum_{\iota=1}^{i-1}\mathcal{E}_{j,\iota}(t)+\sum_{\iota=j+1}^{j+2}\mathcal{E}_{\iota, i-2}(t),\ \ \ i\geq2.\label{q1est}
\end{align}
For $Q_2,$ it follows from \eqref{apriori} \eqref{j1} \eqref{j2} and \eqref{sgy} that
\begin{align}
|Q_2|\les& \epsilon_0\sum_{\iota=0}^{i-1}\Big(|\pt^j\py^{i-1-\iota}(y\sigma \zeta_{yy})|+|\pt^j\py^{i-1-\iota}(\sigma \zeta_y)|\notag\\
&\ \ \ \ \ \ \ \ \ +|\pt^j\py^{i-1-\iota}(y\sigma_y \zeta_y)|+|\pt^{j}\py^{i-1-\iota}(\sigma_y\zeta)|\Big)\notag\\
\les& \epsilon_0\sum_{\iota=0}^{i-1}\Big(|y\sigma \pt^j\py^{i+1-\iota}\zeta|+\sum_{m=0}^{i-\iota}|\pt^j \py^m \zeta | \Big),
\end{align}
which implies that 
\begin{align}
\norm{y \sigma^{\frac{\alpha+i-1}{2}} Q_2}^2 \les \epsilon^2_0\Big(\mathcal{E}_{j,i}(t)+\sum_{\iota=0}^{i-1}
\mathcal{E}_{j, \iota}\Big). \label{q2est}
\end{align}
Here the fact   $(n-[n/2]-\alpha)/2>1/2$ for $\ga>4/3$ have been used.
So, for suitably small $\epsilon_0$, we deduce from \eqref{ejiq12} \eqref{q1est} \eqref{q2est} that 
\begin{align}
\mathcal{E}_{j,i}(t)\les \mathcal{E}_j(t)+\mathcal{E}_{j+1}(t), \ \ \ i=1,\label{ejiej1}
\end{align}
and
\begin{align}
\mathcal{E}_{j,i}(t)\les  \mathcal{E}_{j, i-1}(t)+ \mathcal{E}_{j+2, i-2}(t)+\mathcal{E}_{j+1, i-2}(t)+\sum_{\iota=0}^{i-1}
\mathcal{E}_{j, \iota}(t),\ \ \ i\geq 2.\label{ejiej2}
\end{align}

With \eqref{ejiej1} and \eqref{ejiej2} in hand, together with the induction hypothesis \eqref{hypo}, we are ready to 
show \eqref{ellmeji} holds for $j+i=k+1.$  We first choose $j=k$ and $i=1.$ Then it follows from \eqref{ejiej1} that
\begin{align}
\mathcal{E}_{k, 1}(t)\les \mathcal{E}_{k}(t)+\mathcal{E}_{k+1}(t).\label{induc1}
\end{align}
If we choose $j=k-1$ and $i=2$ in \eqref{ejiej2}, with the aid of \eqref{mej0} and \eqref{hypo},  we obtain 
\begin{align}
\mathcal{E}_{k-1,2}(t)\les  \mathcal{E}_{k-1, 1}(t)+ \mathcal{E}_{k+1, 0}(t)+\mathcal{E}_{k, 0}(t)+
\mathcal{E}_{k-1, 0}(t)\les \sum_{\iota=0}^{k+1} \mathcal{E}_{\iota}(t).\label{induc2}
\end{align}
For other cases, $j=k-2$ and $i=3,$ $j=k-3$ and $i=4$, and other pairs $(j, i)$ satisfying $j+i=k+1,$  with the help of \eqref{induc1} and \eqref{induc2}, we can handle them  similarly.  Now we finish the proof of Lemma \ref{le2}.

\end{proof}

\subsection{Basic Energy Estimate}
\begin{lemma}\label{lemma-1}
Assume that \eqref{apriori} holds for suitably small $\epsilon_0\in (0,1).$ If $\ga>4/3,$ then it holds
\begin{align}\label{le-1}
e^{\delta t}\mathcal{E}_0(t)+ \int_0^te^{\delta s}\mathcal{D}_0(s) ds\leq \mathcal{E}_0(0),\ \ \ \forall t\in[0, T],
\end{align}
for some suitably small $\delta>0.$
\end{lemma}
\begin{proof}
We integrate the product of \eqref{erreqn} and $y^3 \zeta_t$ with respective to the space
variable, by virtue of integration by parts,
 to obtain
\begin{align}
\frac{d}{dt}\int \frac{1}{2}y^4\bar\rho \zeta_t^2 dy
+\int y^4 \bar\rho \zeta_t^2 dy+\int \bar\rho^\ga H_1 dy=0,\label{zt-1}
\end{align}
where
\begin{align*}
H_1=-(1+\zeta)^{-2\ga}(1+\zeta+y\zeta_y)^{-\ga}\lt((1+\zeta)^2y^3\zeta_t\rt)_y 
+\lt((1+\zeta)^{-2}y^3 \zeta_t\rt)_y.
\end{align*}

Indeed,  using a direct calculation,  we  have  the following equality: 
\begin{align}
&\int \bar\rho^{\ga}H_1dy
=\frac{d}{dt}\int y^2 \bar\rho^\ga {E}_0 dy,
\end{align}
where 
\begin{align*}
{E}_0(y,t)=\frac{1}{\ga-1}[&(1+\zeta)^{2-2\ga}(1+\zeta+y\zeta_y)^{1-\ga}\\
&-3(\ga-1)(1+\zeta)^{-1}+
(\ga-1)(1+\zeta)^{-2}y\zeta_y+(3\ga-4)].
\end{align*}
Then we can rewrite  \eqref{zt-1} as following 
\begin{align}
\frac{d}{dt}\int \lt[\frac{1}{2}y^4\bar\rho \zeta_t^2 +y^2\bar\rho^\ga {E}_0(y, t)\rt] dy
+\int y^4 \bar\rho \zeta_t^2 dy=0.\label{zte}
\end{align}
Taking advantage of the smallness of $\zeta$ and $y \zeta_y,$  
it follows from Taylor expansion that
\begin{align}
{E}_0=&\lt[(\frac{9}{2}\ga-6)\zeta^2+(3\ga-4)\zeta y\zeta_y
+\frac{\ga}{2}y\zeta_y^2 \rt]\notag\\
&+O(1)(|\zeta|+|y\zeta_y|)(\zeta^2+y^2\zeta_y^2).\notag
\end{align}

Hence, if $\ga>4/3,$  we have 
\begin{align}
{E}_0\sim&\lt( \zeta^2+(y\zeta_y)^2\label{zt-2}\rt)
\end{align} 
for suitably  small $\zeta$ and $y\zeta_y.$

Multiplying \eqref{erreqn} by $y^3 \zeta$, and then integrating the product with respect to 
the space variable, we apply integration by parts to obtain
\begin{align}
\frac{1}{2}\frac{d}{dt} \int  (2 y^4 \bar\rho \zeta_t\zeta+y^4 \bar\rho \zeta^2) dy +\int \bar\rho^\ga H_2 dy
=\int y^4 \bar\rho \zeta_t^2 dy,\label{z-1}
\end{align}
where 
$$H_2:=(1+\zeta)^{-2\ga}(1+\zeta+y\zeta_y)^{-\ga}((1+\zeta)^2y^3 \zeta)_y-((1+\zeta)^{-2}y^3 \zeta)_y,$$
due to $\zeta(r_0, t)=0, \sigma(R)=0.$
It follows from Taylor expansion and the smallness of $\zeta$ and $y \zeta_y$ that
\begin{align}
H_2\geq& y^2 [(9\ga-12)\zeta^2+(6\ga-8)y\zeta_y \zeta+\ga(y\zeta_y)^2]
-C\epsilon_0y^2(\zeta^2+y^2\zeta_y^2)\notag\\
\ges & y^2 (\zeta^2+y^2\zeta_y^2), \label{h2}
\end{align}
where the second inequality follows from $\ga>4/3$ and the smallness of  $\epsilon_0.$

Calculate $\eqref{z-1}+2\times\eqref{zte}$ to obtain
\begin{align}
\frac{d}{dt} \mathfrak{E}_0(t)+\mathfrak{D}_{0}(t)=0,\label{ded}
\end{align}
where
\begin{align}
&\mathfrak{E}_0(t)= \int  (y^4 \bar\rho \zeta_t\zeta+\frac{1}{2} y^4 \bar\rho \zeta^2) dy+\int \lt[y^4\bar\rho \zeta_t^2 +2y^2\bar\rho^\ga {E}_0(y, t)\rt] dy,\notag\\
&\mathfrak{D}_{0}(t)=\int (y^4 \bar\rho\zeta_t^2 +\bar\rho^\ga H_2)dy.\notag
\end{align}
Clearly, if $\ga>4/3$ and $\epsilon_0$ is suitably small, it follows from \eqref{zt-2} and \eqref{h2} that
\begin{align}
\mathfrak{E}_0\sim \mathcal{E}_0,\ \ \mathfrak{D}_0\sim \mathcal{D}_0,\label{eedd}
\end{align}
and 
\begin{align}
\mathcal{E}_0(t)\les \mathcal{D}_0(t),\label{e0d0}
\end{align}
due to 
\begin{align*}
\int y^4 \sigma^\alpha \zeta^2 dy\les& \int_{r_0}^{(r_0+R)/2} y^2 \sigma^{\alpha+1} \zeta^2 dy
+\int_{(r_0+R)/2}^R \sigma^{\alpha} \zeta^2 dy\\
\les& \mathcal{D}_0+ \int_{(r_0+R)/2}^R \sigma^{\alpha+2} (\zeta^2 +\zeta_y^2)dy\\
\les & \mathcal{D}_0+ \int_{(r_0+R)/2}^Ry^2 \sigma^{\alpha+2} \lt(\zeta^2 + (y\zeta_y)^2\rt)dy\\
\les & \mathcal{D}_0.
\end{align*}

Then, with the aid of \eqref{eedd} and \eqref{e0d0}, we integrate the product of $e^{\delta t}$ and \eqref{ded} with respect to time variable to obtain \eqref{le-1} for suitably small $\delta$. 
Hence  we complete the proof of  Lemma \ref{lemma-1} .

%\begin{align}
%&\int [y^4 \bar\rho (\zeta^2+\zeta_t^2)+y^2 \bar\rho^\ga ( \zeta^2+y^2 \zeta_y^2)](y,t)dy \notag\\
%&+\int_0^t \int \lt( y^4\bar\rho \zeta_s^2+y^2\bar\rho^\ga (\zeta^2+y^2 \zeta_y^2) \rt)(y,s)dyds
%\les \mathcal{E}^{0,0}(0),
%\end{align}
%where \eqref{zt-2} and \eqref{h2} have been used. 

\end{proof}

\subsection{Higher-Order Energy Estimates}   
The  equation \eqref{erreqn} can be wrote  as 
\begin{align}
&y \bar\rho\zeta_{tt}+y\bar\rho \zeta_t+\lt\{\bar\rho^\ga[(1+\zeta)^{2-2\ga}(1+\zeta+y\zeta_y)^{-\ga}-(1+\zeta)^{-2}]\rt\}_y\notag\\
&-2\bar\rho^\ga[(1+\zeta)^{1-2\ga}(1+\zeta+y\zeta_y)^{-\ga}+(1+\zeta)^{-3}]\zeta_y=0.\label{reerr}
\end{align}
The linearized  equation of \eqref{reerr} reads 
\begin{align*}
y \bar\rho\zeta_{tt}+y\bar\rho \zeta_t+\lt\{\bar\rho^\ga[(4-3\ga)\zeta-\ga y\zeta_y]\rt\}_y-4\bar\rho^\ga \zeta_y=0.
\end{align*}
Let $j\geq 1$ and take the $j$-th time derivatives of the equation \eqref{reerr} to obtain 
\begin{align}
y \bar\rho&\pt^j\zeta_{tt}+y\bar\rho \pt^j\zeta_t+\lt[\bar\rho^\ga \lt(h_1\pt^j\zeta+h_2y\pt^j\zeta_y+R_1\rt)\rt]_y\notag\\
&+\bar\rho^\ga\lt[(3h_2-h_1)\pt^j\zeta_y+R_2\rt]-2\bar\rho^\ga\lt(h_3\zeta_y\pt^j\zeta+R_3\rt)=0,\label{ptj}
\end{align}
where
\begin{align}
h_1=&(2-2\ga)(1+\zeta)^{1-2\ga}(1+\zeta+y\zeta_y)^{-\ga}\notag\\
&-\ga(1+\zeta)^{2-2\ga}(1+\zeta+y\zeta_y)^{-\ga-1}+2(1+\zeta)^{-3},\notag\\
h_2=&-\ga (1+\zeta)^{2-2\ga}(1+\zeta+y\zeta_y)^{-\ga-1},\notag\\
h_3=&(1-2\ga)(1+\zeta)^{-2\ga}(1+\zeta+y\zeta_y)^{-\ga}\notag\\
&-\ga(1+\zeta)^{1-2\ga}(1+\zeta+y\zeta_y)^{-\ga-1}-3(1+\zeta)^{-4},\notag
\end{align}
and 
\begin{align*}
R_1=&\pt^{j-1}(h_1 \zeta_t+h_2y\zeta_{yt})-\lt(h_1\pt^j\zeta+h_2 y\zeta\pt^j\zeta_y\rt),\\
R_2=&\pt^{j-1}[(3h_2-h_1)\zeta_{ty}]-(3h_2-h_1)\pt^j \zeta_y,\\
R_3=&\pt^{j-1}(h_3\zeta_y\zeta_t)-h_3\zeta_y\pt^j\zeta.
\end{align*}
Here $R_1, R_2$ and $R_3$ refer to lower-order terms involving $\pt^\iota (\zeta, y\zeta_y)$ with 
$\iota=0, 1, \cdots, j-1.$ Taking advantage of the smallness of $\zeta$
 and $y\zeta_y$, one can apply Taylor expansion on $h_\iota,\iota=1,2,3,$  to obtain
 \begin{align}
 h_1=&(4-3\ga)+(9\ga^2-9\ga-4)\zeta+\ga(3\ga-1)y\zeta_y+\bar h_1,\notag\\
 h_2=&-\ga+\ga(3\ga-1)\zeta+\ga(\ga+1)y\zeta_y+\bar h_2,\notag\\
 h_3=&-2-3\ga+\bar h_3, \label{hhbar}
 \end{align}
 where $\bar h_\iota, \iota=1,2,3,$ satisfy
 \begin{align}
 |\bar h_1|+|\bar h_2|\les\lt(\zeta^2+(y\zeta_y)^2\rt),\ \ \ |\bar h_3|\les (|\zeta|+|y\zeta_y|).\label{hbarest}
 \end{align}
 
 \begin{lemma}\label{lemma-2}
 Assume that \eqref{apriori} holds for suitably small $\epsilon_0\in (0,1).$ If $\ga>4/3,$ then, for $1\leq j\leq n,$ it holds
 \begin{align}
 e^{\delta t}\mathcal{E}_j(t)+\int_0^t e^{\delta s}\mathcal{D}_j(s)ds\les \mathcal{E}(0),\ \ \ \forall t\in[0, T]. \label{le2-1}
 \end{align}
 \end{lemma}
 
 \begin{proof}
 We integrate the product of \eqref{ptj} and $y^3\pt^j\zeta_t$ with respect to the spatial variable, with the aid of $\pt^{j}\zeta(r_0, t)=0,$  to obtain
 \begin{align}
 \frac{d}{dt}\lt[\int \lt(\frac{1}{2}y^4 \bar\rho \lt(\pt^j \zeta_t\rt)^2 +y^2 \bar\rho^\ga E_j \rt)dy+ \tilde{R}^{j-1}\rt]+\int y^4 \bar\rho\lt(\pt^j \zeta_t\rt)^2 dy=F_1^j+F_2^j,\label{ptjzt}
 \end{align}
 where 
 \begin{align}
 &E_j=-\frac{1}{2}\lt[(3h_1+2h_3 y\zeta_y)\lt(\pt^j\zeta\rt)^2+2h_1 (\pt^j \zeta)\lt(y\pt^j\zeta_y\rt)+h_2\lt(y\pt^j \zeta_y\rt)^2 \rt],\notag\\
 &\tilde{R}^{j-1}=-\int y^2 \bar\rho^\ga \lt[(3R_1\pt^j \zeta+R_1 y\pt^j \zeta_y)-y(R_2-2R_3)\pt^j\zeta\rt]dy,\notag\\
 &F_1^j=-\int \frac{1}{2} y^2 \bar\rho^\ga \lt[(3h_1+2h_3y\zeta_y)_t(\pt^j \zeta)^2+
 2h_{1t}\lt(y\pt^j\zeta_y\rt)\pt^j\zeta+h_{2t}\lt(y\pt^j\zeta_y\rt)^2\rt] dy,\notag\\
 &F_2^j=-\int \bar\rho^\ga y^2\lt[3R_{1t}\pt^j \zeta+R_{1t}y\pt^j\zeta_y-(R_2-2R_3)_t\pt^j\zeta\rt]dy.\notag
 \end{align}
 In particular, $\tilde{R}^{j-1}=0$ when $j=1,$ since $R_1=R_2=R_3=0$ when $j=1.$

 Similarly, we integrate the product of \eqref{ptj} and $y^3 \pt^j\zeta$ with respect to spatial variable to get 
 \begin{align}
 \frac{d}{dt}\int y^4 \bar\rho \lt( \pt^j\zeta \pt^j\zeta_t+\frac{1}{2}\lt(\pt^j\zeta\rt)^2\rt)dy -\int y^4\bar\rho\lt(\pt^j\zeta_t\rt)^2dy +\int y^2 \bar\rho^\ga D_j dy+\tilde{R}^{j-1}=0,\label{ptjz}
 \end{align}
 where 
 \begin{align}
 D_j=&-(h_1\pt^j\zeta +h_2y\pt^j\zeta_y)(3\pt^j \zeta+y\pt^j\zeta_y)\notag\\
 &+(3h_2-h_1)\lt(y\pt^j\zeta_y\rt)\pt^j \zeta
 -2h_3 (y\zeta_y)\lt(\pt^j \zeta\rt)^2.\label{dj}
 \end{align}
 Here $D_j$  satisfies that 
 \begin{align}
 D_j\geq &3(3\ga-4)(\pt^j\zeta)^2+2(3\ga-4)\pt^j\zeta\lt(y\pt^j\zeta_y\rt)+\ga \lt(y\pt^j\zeta_y\rt)^2\notag\\
 &-C(|\zeta|+|y\zeta_y|)\lt[\lt(\pt^j\zeta\rt)^2+\lt(y\pt^j\zeta_y\rt)^2\rt]\notag\\
 \ges & \lt[\lt(\pt^j\zeta\rt)^2+\lt(y\pt^j\zeta_y\rt)^2\rt],\ \ \
 {\rm if} \ \ \ga>\frac{4}{3}, \label{djest1}
 \end{align}
 and
 \begin{align}
 D_j\les& \lt[\lt(\pt^j\zeta\rt)^2+\lt(y\pt^j\zeta_y\rt)^2\rt],\label{djest2}
 \end{align}
 due to \eqref{hhbar} \eqref{hbarest}  \eqref{apriori}and suitably small $\epsilon_0.$
 
 Then, we calculate $\eqref{ptjz}+2\times \eqref{ptjzt}$ to obtain 
 \begin{align}
 \frac{d}{dt}\lt(\mathfrak{E}_j+2\tilde{R}^{j-1}\rt)+\mathfrak{D}_j=2(F_1^j+F_2^j)-\tilde{R}^{j-1},\label{energysum}
 \end{align}
 where 
 \begin{align}
 \mathfrak{E}_j=&\int \lt\{ y^4 \bar\rho \lt[\lt(\pt^j \zeta_t\rt)^2+ \pt^j\zeta \pt^j\zeta_t+\frac{1}{2}\lt(\pt^j\zeta\rt)^2\rt] +2y^2 \bar\rho^\ga E_j \rt\}dy,\label{mfej}\\
 \mathfrak{D}_j=&\int y^4\bar\rho\lt(\pt^j\zeta_t\rt)^2dy +\int y^2 \bar\rho^\ga D^j dy.\label{mfdj}
 \end{align}
  It follows from Taylor expansion and the smallness of $\zeta$ and $y\zeta_y$ that
 \begin{align}
 E_j=&\frac{1}{2}\lt[3(3\ga-4)\lt(\pt^j\zeta\rt)^2+2(3\ga-4)\pt^j\zeta\lt(y\pt^j \zeta_y\rt)+\ga\lt(y\pt^j\zeta_y\rt)^2\rt]\notag\\
 &+O(1)(|\zeta|+|y\zeta_y|)\lt((\pt^j\zeta)^2+(y\pt^j \zeta_y)^2\rt)\notag\\
 \sim &\lt((\pt^j\zeta)^2+(y\pt^j \zeta_y)^2\rt).\label{ejsim}
 \end{align}
 This, together with \eqref{djest1} \eqref{djest2} \eqref{mfej} and \eqref{mfdj}, implies that 
 \begin{align}
 \mathfrak{E}_j(t)\sim \mathcal{E}_j(t),\ \ \ {\rm and}\ \ \ \mathfrak{D}_j(t)\sim\mathcal{D}_j(t), \label{eeddj}
 \end{align}
 if $\ga>4/3$ and $\epsilon_0$ is suitably small.
 
 Furthermore, similar to derive \eqref{e0d0}, it follows from  
 the definition of $\mathcal{E}_j,$  $\mathcal{D}_j$ and $\mathcal{D}$ that 
\begin{align}
&\mathcal{E}_j(t)\les \mathcal{D}_j(t),\label{ejdj}\\
&\sum_{j=0}^n \mathcal{E}_j\les  \int y^4 \sigma^\alpha \zeta^2 dy+\mathcal{D}(t)\les \mathcal{D}(t).\label{sumejd}
\end{align}
Thus,  we can deduce from \eqref{sumejd} and \eqref{proest} and the definition of $\mathcal{E}(t)$ that
\begin{align}
\mathcal{E}(t)\les \mathcal{D}(t).\label{eleqd}
\end{align}

 Notice that
 \begin{align}
|F_1^j|\les& \int y^2 \bar\rho^\ga \lt[(|\zeta|+|y\zeta_y|)+(|\zeta|+|y\zeta_y|)(|\zeta_t|+|y\zeta_{yt}|)\rt]\lt(\lt(\pt^j\zeta\rt)^2+\lt(y\pt^j\zeta_y\rt)^2\rt)dy\notag\\
\les&\epsilon_0 \int y^2 \bar\rho^\ga \lt(\lt(\pt^j\zeta\rt)^2+\lt(y\pt^j\zeta_y\rt)^2\rt)dy\les \epsilon_0 \mathcal{E}_j,\label{f1j}
 \end{align}
 where Taylor expansion and \eqref{apriori} have been used.
 
 {\it Claim}: for $F_2^j$ and $\tilde{R}^{j-1}, j\geq 1,$  it holds that 
 \begin{align}
  & |F_2^j|\les \epsilon_0 \mathcal{E}_j+ \epsilon_0\sum_{\iota=0}^{j-1} \mathcal{E}_{\iota},\label{f2jest}\\
 & |\tilde{R}^{j-1}|\les \epsilon_0 \mathcal{E}_j+ \epsilon_0 \sum_{\iota=0}^{j-1} \mathcal{E}_{\iota},\label{trjest}
 \end{align}
where the a priori estimate \eqref{apriori} has been used.  We will prove this claim in the Appendix.
 
 Once we have \eqref{f2jest} and \eqref{trjest}, together with \eqref{f1j},  we multiply \eqref{energysum} by $e^{\delta t}$ to obtain
  \begin{align}
  &\frac{d}{dt}\lt(e^{\delta t} (\mathfrak{E}_j+2\tilde{R}^{j-1})\rt)-\delta e^{\delta t} (\mathfrak{E}_j+2\tilde{R}^{j-1})+e^{\delta t}\mathfrak{D}_j\notag\\
  &\les \epsilon_0  e^{\delta t} \mathcal{E}_j+\epsilon_0 e^{\delta t}\sum_{\iota=0}^{j-1}\mathcal{E}_\iota.\label{edttotal}
  \end{align}
  When $j=1,$ one has 
  \begin{align}
  &\frac{d}{dt}\lt(e^{\delta t} \mathfrak{E}_1\rt)-\delta e^{\delta t} \mathfrak{E}_1+e^{\delta t}\mathfrak{D}_1
  \les \epsilon_0  e^{\delta t} \mathcal{E}_1+\epsilon_0 e^{\delta t}\mathcal{E}_0,\label{edt1}
 \end{align}
 which,  together with \eqref{le-1} \eqref{e0d0} \eqref{eeddj}, implies 
 \begin{align}
 e^{\delta t} \mathcal{E}_1 +\int_0^t e^{\delta s} \mathcal{D}_1ds \les \mathcal{E}_1(0),\label{edt2}
 \end{align}
 for suitably small $\delta$ and $\epsilon_0$.
 
 When $j=2,$ it follows from \eqref{eeddj} \eqref{ejdj} \eqref{edttotal} \eqref{f2jest} and \eqref{trjest} that, for suitably small 
 $\delta$ and $\epsilon_0$,
 \begin{align}
  &e^{\delta t} \mathcal{E}_2(t)+\int_0^t e^{\delta s}\mathcal{D}_2(s) ds\notag\\
  &\les  e^{\delta t}\sum_{\iota=0}^{1}\mathcal{E}_\iota(t)+\int_0^s \lt( e^{\delta s}\sum_{\iota=0}^{1}\mathcal{E}_\iota(s) \rt)ds\les \mathcal{E}(0),\label{edt3}
  \end{align}
where \eqref{le-1} \eqref{ejdj} and \eqref{edt2} have been used for the second inequality in \eqref{edt3}.
Similarly, for other cases $k\geq 2,$ using bootstrap argument,    one has
 \begin{align}
  &e^{\delta t} \mathcal{E}_k(t)+\int_0^t e^{\delta s}\mathcal{D}_k(s) ds\notag\\
  &\les  e^{\delta t}\sum_{\iota=0}^{k-1}\mathcal{E}_\iota(t)+\int_0^s \lt( e^{\delta s}\sum_{\iota=0}^{k-1}\mathcal{E}_\iota(s) \rt)ds\les \mathcal{E}(0),\label{edt4}
  \end{align}
  for suitably small $\delta$ and $\epsilon_0$.
  Thus we can conclude \eqref{le2-1} by virtue of \eqref{le-1} and the definition of $\mathcal{D}_j(t)$ and 
  $\mathcal{E}_j(t).$ Now the proof of Lemma \ref{lemma-2} is completed.
  
 \end{proof}
 {\textit {Proof of Theorem \ref{thm1} and Theorem \ref{thm2}.}} It follows from \eqref{le-1} \eqref{le2-1} and  the elliptic estimates \eqref{proest} and \eqref{sumejd} \eqref{eleqd} that 
  \begin{align}
 e^{\delta t}\mathcal{E}(t)+\int_0^t e^{\delta s}\mathcal{D}(s)ds\les \mathcal{E}(0),\ \ \ \forall t\in[0, T], \label{thm1p}
 \end{align}
 if $\ga>4/3$ and $r_0< R\leq \frac{4}{3-\alpha}r_0<\infty.$ Then we obtain \eqref{thm1-1} and finish the proof of Theorem \ref{thm1} by a standard global existence argument.
 
 Note that 
 \begin{align*}
 \rho(t, \eta(t,y))-\bar\rho(y)=&\frac{\bar \rho(y)}{(1+\zeta)^2(1+\zeta+y\zeta_y)}-\bar\rho(y)\\
 =&-\bar\rho(y)\frac{(\zeta+y\zeta_y)+\zeta(2+\zeta)(1+\zeta+y\zeta_y)}{(1+\zeta)^2(1+\zeta+y\zeta_y)},\\
 u(t, \eta(t, y))=\eta_t(t,y)=&y\zeta_{t},\ \ \ R(t)-R=R\zeta(R, t).
 \end{align*}
 This, together with \eqref{thm1p}  \eqref{zeng2017} and the fact $H^{1/2+\varepsilon}(I)\hookrightarrow L^{\infty}(I)$ for $\varepsilon>0$,  implies \eqref{rhodec} and \eqref{urdec}. Then the proof of
 Theorem \ref{thm2} is completed.
 
  {\textit {Proof of Theorem \ref{thm3}.}}
  Here we note that the equilibria equation \eqref{gravity1} has a unique solution only if $\ga\geq 4/3.$ There are  multiple solutions to equation \eqref{gravity1} with $1<\ga<4/3,$ which  has been shown in \cite{KL}. The uniqueness of the equilibrium for Euler-Poisson equations with damping and solid core  requires $\ga\geq 4/3,$  compared with that $\ga>1$ for Euler equations with damping and solid core.
  
Similar to the strategy in the proof of Theorem \ref{thm1} and Theorem \ref{thm2}, by means of elliptic estimates,  weighted temporal energy estimates and weighted Sobolev embedding theory, one can prove the  global existence of smooth solution and stability property to  Euler-Poisson equations with damping and solid core  \eqref{epc}  for  any given  $\ga>4/3$ provided that the initial spherically symmetric perturbation around the equilibrium \eqref{gravity1} is sufficiently small.  We omit the repeat for simplicity.

\subsection*{Acknowledgements}
This research was supported in part by China Postdoctoral  Science Foundation  grant 2021M691818;  Natural Science Foundation of China  (NSFC)   Grants 12101350, 12171267, 12271284.
The author is deeply grateful to the referee for his/her  helpful suggestions which have helped to clarify some important points and improved
the quality  of the paper greatly.

\subsection*{Data Availability} 
Data sharing not applicable to this article as no datasets were generated or analysed during the current study.

 \vspace{1cm}

\section*{Appendix}

\noindent{\it Claim}:   For $F_2^j$ and $\tilde{R}^{j-1}, j\geq 1,$ it holds that 
 \begin{align}
 & |F_2^j|\les \epsilon_0 \mathcal{E}_j+ \epsilon_0\sum_{\iota=0}^{j-1} \mathcal{E}_{\iota},\label{af2jest}\\
 & |\tilde{R}^{j-1}|\les \epsilon_0 \mathcal{E}_j+ \epsilon_0 \sum_{\iota=0}^{j-1} \mathcal{E}_{\iota},\label{atrjest}
 \end{align}
 for a suitably small positive constant $\epsilon_0.$
 
\noindent{\it Proof of the Claim.}
If $j=1,$ then $R_1=R_2=R_3=\tilde{R}^0=0,$ it is trivial.

For $F_2^j,  6\geq n\geq j\geq 2$ (since $\alpha<3$ for $\ga>4/3$),  note that
\begin{align}
&R_{1t}=(j-1)\lt(h_{1t}\pt^j \zeta+h_{2t}y\pt^j\zeta_{y}\rt)+O(1)\sum_{\iota=2}^{j}|\pt^{\iota}(h_1, h_2)|\lt|\pt^{j+1-\iota}(\zeta, y\zeta_y)\rt|,\notag\\
&yR_{2t}=(j-1)(2h_2-h_1)_t\lt(y\pt^j\zeta_y\rt)+O(1)\sum_{\iota=2}^{j}|\pt^{\iota}(h_1, h_2)|\lt|\pt^{j+1-\iota}(y\zeta_y)\rt|,\notag\\
&yR_{3t}=(j-1)(h_3 y\zeta_y)_t\lt(\pt^j\zeta\rt)+O(1)\sum_{\iota=2}^{j}|\pt^{\iota}(h_3y\zeta_y)|
\lt|\pt^{j+1-\iota}\zeta\rt|.\notag
\end{align}
This, together with 
\eqref{hhbar} \eqref{apriori}, implies that
\begin{align}
|F_2^j|\les \epsilon_0 \int y^2 \sigma^{\alpha+1} \lt|\pt^j(\zeta, y\zeta_y)\rt|^2dy +\int y^2 \sigma^{\alpha+1}\lt(\lt|\pt^j\zeta\rt|+
\lt|y\pt^j\zeta_y\rt|\rt) \mathfrak{R} dy,\label{f2jR}
\end{align}
where
\begin{align}\mathfrak{R}=\sum_{\iota=2}^{j}\lt|\pt^{\iota}(h_1, h_2, h_3 y\zeta_y)\rt|\lt|\pt^{j+1-\iota}(\zeta, y\zeta_y)\rt|.\label{frkR}
\end{align}
Here $\mathfrak{R}$ satisfies that
\begin{align}
\mathfrak{R}\les \epsilon_0 \sum_{\iota=1}^j\lt|\pt^\iota (\zeta, y\zeta_y)\rt|+\overline{\mathfrak{R}},
\end{align}
where 
\begin{align}
\overline{\mathfrak{R}}=&\lt|\pt^2(\zeta, y\zeta_y)\rt|\lt|\pt^{j-1}(\zeta, y\zeta_y)\rt|+\lt(\lt|\pt^2(\zeta, y\zeta_y)\rt|+\lt|\pt^3(\zeta, y\zeta_y)\rt| \rt)\lt|\pt^{j-2}(\zeta, y\zeta_y)\rt|\notag\\
&+\lt(\lt|\pt^2(\zeta, y\zeta_y)\rt|+\lt|\pt^3(\zeta, y\zeta_y)\rt|+\lt|\pt^4(\zeta, y\zeta_y)\rt| \rt)\lt|\pt^{j-3}(\zeta, y\zeta_y)\rt|,\label{bfrkR}
\end{align}
since $j\leq n\leq 6.$
Then we use \eqref{apriori} to obtain 
\begin{align}
\overline{\mathfrak{R}}\les \epsilon_0 \sum_{\iota=1}^{[(j-1)/2]} \sigma^{-\frac{\iota}{2}}\lt|\pt^{j-\iota}(\zeta, y\zeta_y)\rt|.\label{bfrkRest}
\end{align}
Now we have 
\begin{align}
&\int y^2 \sigma^{\alpha+1}\lt(\lt|\pt^j\zeta\rt|+\lt|y\pt^j\zeta_y\rt|\rt) \overline{\mathfrak{R}} dy\notag\\
&\les \epsilon_0\int y^2 \bar\rho^\ga \lt|\pt^j(\zeta, y\zeta_y)\rt|^2dy+\epsilon_0\sum_{\iota=1}^{[(j-1)/2]}\int y^2  \sigma^{\alpha+1-\iota}\lt|\pt^{j-\iota}(\zeta, y\zeta_y)\rt|^2dy.\label{Rest1}
\end{align}
For $\iota=1, \cdots, [(j-1)/2],$  it follows from Hardy inequality \eqref{hardy} and the fact $\alpha+1-\iota\geq \alpha+1-[(1+[\alpha])/2]>0$ that 
\begin{align}
&\int y^2  \sigma^{\alpha+1-\iota}\lt|\pt^{j-\iota}(\zeta, y\zeta_y)\rt|^2dy\notag\\
&\les \int_{I_l} y^2 \sigma^{\alpha+1}\lt|\pt^{j-\iota}(\zeta, y\zeta_y)\rt|^2dy+\int_{I_b}  \sigma^{\alpha+1-\iota}\lt|\pt^{j-\iota}(\zeta, y\zeta_y)\rt|^2dy\notag\\
&\les \mathcal{E}_{j-\iota}(t)+\int_{I_b} \sigma^{\alpha+3-\iota}\lt|\pt^{j-\iota}(\zeta, \zeta_y, \zeta_{yy})\rt|^2dy\notag\\
&\les\mathcal{E}_{j-\iota}(t)+\cdots\notag\\
&\les\mathcal{E}_{j-\iota}(t)+\sum_{i=0}^{\iota+1}\int_{I_b} y^4 \sigma^{\alpha+1+\iota}\lt|\pt^{j-\iota}\py^i \zeta\rt|^2 dy
\notag\\
&\les \mathcal{E}_{j-\iota}(t)+\sum_{i=1}^{\iota}\mathcal{E}_{j-\iota, i}(t)\notag\\
&\les \sum_{\iota=0}^{j} \mathcal{E}_{\iota}(t),\label{Rest2}
\end{align}
where the elliptic estimate \eqref{proest} has been used in the last inequality. Then \eqref{af2jest} follows from \eqref{f2jR}-\eqref{Rest2}. 

As for $\tilde{R}^{j-1}, j\geq 2,$ it can be handled by a similar but much easier way to dealing with $F_2^j.$ Hence we omit the repeat for simplicity. Now we finish the proof of the claim stated in \eqref{f2jest} and \eqref{trjest}.

\end{document}